\documentclass[a4paper]{article}

\usepackage[a4paper,top=3cm,bottom=2cm,left=3cm,right=3cm,marginparwidth=1.75cm]{geometry}

\usepackage[english]{babel}
\usepackage[T1]{fontenc} %
\usepackage[utf8x]{inputenc}
\usepackage{amsmath}
\usepackage{amsthm}
\usepackage{thmbox}
\usepackage{stmaryrd}
\usepackage{graphicx}
\usepackage{amssymb}
\usepackage{float}
\usepackage[shortlabels]{enumitem}
\usepackage[colorinlistoftodos]{todonotes}
\usepackage[colorlinks=true, allcolors=black]{hyperref}
\usepackage[export]{adjustbox}
\usepackage[
backend=biber,
style=alphabetic,
sorting=ynt
]{biblatex}
\usepackage{pgf,tikz}
\usetikzlibrary{arrows}

\newtheorem{thm}{Theorem}[section]
\newtheorem{lem}[thm]{Lemma}
\newtheorem{mydef}[thm]{Definition}
\newtheorem{cor}[thm]{Corollary}
\newtheorem{prop}[thm]{Proposition}
\newtheorem{cl}{Claim}
\newtheorem{rem}[thm]{Remark}

\newcommand{\llangle}[0]{\langle\hspace{-.15em}\langle}

\newcommand{\rrangle}[0]{\rangle\hspace{-.15em}\rangle}

\title{Linearity of fundamental groups of graphs of virtually cyclic groups}
\author{Hsuan-Yu Wang}
\begin{document}

\maketitle
\begin{abstract}

We characterize when a generalized Baumslag-Solitar group is linear, and extend the result to the fundamental groups of a graph of groups with infinite virtually cyclic vertex and edge groups.

\end{abstract}
\section{Introduction}

A generalized Baumslag-Solitar group (GBS group) is the fundamental group of a graph of groups where all vertex groups and edge groups are infinite cyclic. As the name suggests, GBS groups generalize the classical Baumslag–Solitar groups (BS groups), which are a rich source of (counter)examples in geometric group theory. They are given by the presentation

\[BS(m,n)=\langle x,t\mid t^{-1}x^mt = x^n\rangle.\]

These groups were defined by Gilbert Baumslag and Donald Solitar in 1962 to provide examples of non-Hopfian groups. However, some Baumslag–Solitar groups are Hopfian, and so it is important to characterize the conditions under which a BS group is (non-)Hopfian. It turns out that a BS group is Hopfian if and only if is residually finite if and only if it is linear if and only if $|n|=1$ or $|m|=1$ or $|n|=|m|$.

Our first goal is to characterize when a GBS group is linear (and we will see that it coincides precisely when it is residually finite and Hopfian). More precisely, we prove the following theorem.

\begin{thm}\label{thm:main}
Let $G$ be a generalized Baumslag-Solitar group. Let $\langle x \rangle$ be a vertex group. Then, at least one of the following alternatives holds:
\begin{enumerate}[(a)]
\item for some integer $k\geq 1$, $\langle x^k\rangle $ is a normal subgroup.
\item $G$ is a solvable Baumslag-Solitar group.
\item Otherwise, $G$ is not residually finite.
\end{enumerate}
Furthermore, a GBS group $G$ is virtually a direct product of $\mathbb{Z}$ and a free group if it satisfies condition (a), and it satisfies condition (b) if and only if $G = BS(m,n)$ with $\min(|m|,|n|)=1$. In both of these conditions, $G$ is linear (and hence residually finite).

In addition, a GBS group $G$ satisfies condition (a) if and only if it is unimodular (see \ref{def:unimodular})
\end{thm}

This theorem was conjectured by \textbf{Yves de Cornulier} on the forum \cite{Cornulier} and was proved in \cite[Corollary 7.7]{Levitt}.

In this paper, we provide a self-contained proof of the theorem, and we will extend it to fundamental groups of graph of groups with infinite virtually cyclic vertex and edge groups. In this case, we have the following theorem.

\begin{thm}\label{thm:v.c.}
Let $G$ be the fundamental group of a graph of groups with infinite virtually cyclic vertex and edge groups. Then, at least one of the following alternatives holds:
\begin{enumerate}[(a)]
\item for some $g \in G$ belonging to a vertex group, $\langle g\rangle $ is a normal subgroup.
\item $G$ is an ascending HNN-extension of an infinite virtually cyclic group $H$ relative to $\varphi: H \rightarrow K$.
\item $G$ is not residually finite.
\end{enumerate}
Furthermore, the group $G$ is virtually a direct product of $\mathbb{Z}$ and a free group if it satisfies condition (a). In both conditions (a) and (b), $G$ is linear (and hence residually finite).
\end{thm}

In Section 2, we prove that if for some integer $k\geq 1$, $\langle x^k\rangle $ is a normal subgroup, then the GBS group is virtually a direct product of $\mathbb{Z}$ and a free group. We also show that a group that is virtually a direct product of $\mathbb{Z}$ and a free group is linear.

In Section 3, we prove the result for the Baumslag-Solitar groups. In addition, we also discuss solvable BS groups.

In Section 4, we prove the case where the graph consists of a single vertex and a finite number of self-loops.

In Section 5, we complete the proof of Theorem \ref{thm:main} for GBS groups.

In Section 6, we extend the result to fundamental groups of graph of virtually cyclic groups.

\section{GBSs with non-trivial cyclic normal subgroup.}

We start this section with an observation that we use throughout the paper.

\begin{lem}\label{lem:cyclic_normal_subgroups}
Let $G$ be a group and $\langle x \rangle$ be an infinite normal subgroup of $G$. Then, for each $t \in G$, we have that $txt^{-1} = t^{-1}xt = x^{\pm 1}$. 
\end{lem}
\begin{proof}
Since $\langle x \rangle$ is normal, we have that for every element $t \in G$, $t^{-1}xt, txt^{-1} \in \langle x \rangle$. It follows that $t^{-1}xt = x^{u}, txt^{-1} = x^{v}$ for some $u,v \in \mathbb{Z}$. But then, $x= t (t^{-1} x t) t^{-1} = tx^u t^{-1} = x^{vu}$ and so $vu = 1$. Since $u,v\in \mathbb Z$, we conclude that $u = v = \pm 1$.
\end{proof}

\begin{lem}\label{lem:BS_normal_linear}
Let $G$ be a GBS group and let $\langle x \rangle$ be a vertex group. If $\langle x^k\rangle$ is a normal subgroup of $G$ for some integer $k\ge 1$, then $G$ is virtually a direct product of $\mathbb Z$ and a free group.    
\end{lem}

\begin{proof}
Let $\langle x^k \rangle$ be a normal subgroup of $G$. By Lemma \ref{lem:cyclic_normal_subgroups}, we have for every edge element $t\in G$ we have $tx^kt^{-1} = x^{\pm k}$, so we see that $\langle x^{k} \rangle$ is contained in all the edge groups, and hence can be embedded in all the vertex groups.

We consider the quotient of $G$ by the normal subgroup $\langle x^k \rangle$ and obtain a fundamental group of a finite graph of finite (cyclic) groups. By \cite[Proposition 11, p.~160]{Serre}, this quotient $G/\langle x^k\rangle$ is virtually free.

Let $H$ be a subgroup of $G$ such that $H/\langle x^k\rangle$ is a free group of finite index in $G/\langle x^k\rangle$. We then consider the group homomorphism $f:H \longrightarrow \text{Aut}(\langle x^k \rangle)$ where $f(h):\langle x^k \rangle \to \langle x^k \rangle$ acts by conjugation, i.e. $f(h)(x^k)=(x^k)^{h}=x^{\pm k} $, by Lemma \ref{lem:cyclic_normal_subgroups}, and let $K$ be its kernel.

Then, $K$ is a normal subgroup of $H$ of index at most 2, so it is a subgroup of $G$ with finite index. Furthermore, note that by definition, for each $g\in K$ we have that $(x^k)^g=x^k$. We next show that $K$ is isomorphic to a direct product of $\mathbb{Z}$ and a free group. We consider the quotient $L = K/ \langle x^k\rangle$; it is a subgroup of $H/\langle x^k\rangle$, so it is also a free group. We then have an exact sequence
\[1\longrightarrow \langle x^k\rangle \longrightarrow K \longrightarrow L\longrightarrow 1.\]
We consider a basis of the free group $L$ and the group homomorphism $r: L\rightarrow K$ which sends every element of this basis to a corresponding representative; the image of this basis must be a basis of a free group (otherwise their projections would not be a basis of the free group $L$). So, the sequence is right-splitting; then, with the splitting lemma, we have that $K$ is isomorphic to $\langle x^k\rangle \rtimes r(L)$. With every element in $K$ commutes with $x^k$ by the definition, we have that this semi-direct product is in fact a direct product, so we have $K$ is isomorphic to the direct product of $\mathbb{Z}$ and the free group $r(L)$.
\end{proof}
\begin{lem} \label{lem:v.l.}
A group $G$ that is virtually linear is linear.
\end{lem}
\begin{proof}
If a subgroup $K$ of $G$ of index $r$ is isomorphic to a subgroup of $GL(n, \mathbb{Q})$, then we have that the induced representation of $G$ is isomorphic to a subgroup of $GL(rn, \mathbb{Q})$. Therefore, $G$ is linear.
\end{proof}

We conclude with a simple observation about Theorem \ref{thm:main}.
\begin{cor}\label{cor:(a)}
A GBS group $G$ that satisfies condition (a) is linear. 
\end{cor}
\begin{proof}
With Lemma \ref{lem:BS_normal_linear}, we have that $G$ is virtually linear because free groups and $\mathbb{Z}$ are linear, so is their direct product. Then we conclude by Lemma \ref{lem:v.l.}
\end{proof}
\section{Baumslag-Solitar Group}
In this section, we prove Theorem \ref{thm:main} for Baumslag-Solitar groups. Let $G$ be $BS(m,n)=\langle x,t\mid t^{-1}x^mt = x^n\rangle$. Since $BS(m,n)=BS(-m,-n)$, we assume without loss of generality that $m > 0$. 

We first recall the well-known Britton's Lemma for HNN-extensions; see \cite[p.~181]{Roger}. We formulate it in the form that we use.

\begin{lem}[Britton's Lemma]
For $G^*$ be an HNN-extension of $G$ relative to an isomorphism $\varphi:H\rightarrow K$ and $w = g_0t^{\varepsilon_1}g_1t^{\varepsilon_2}g_2\cdots t^{\varepsilon_l}g_l$ with $g_i \in G$ and $\varepsilon_i =\pm 1$. If $w =1$, we have that one of the following holds:\\
(a) either $l=0$ and $g_0 = 1$; or\\
(b) $l > 0$ and $t^{\varepsilon_i} g_i t^{\varepsilon_{i+1}}$ is equal to $t^{-1}ht, h\in H$ or $tkt^{-1}, k\in K$, for some $i=1, \dots, l-1$.
\end{lem}

\begin{lem}\label{lem:BS_fr}
If $G$ is residually finite, then one of the following holds
\begin{enumerate}[(a)]
\item  either for some integer $k\geq 1$, $\langle x^k\rangle $ is a normal subgroup, or
\item $G$ is a solvable BS group, i.e. $G=BS(1,n)$.
\end{enumerate}
\end{lem}

Our proof will follow from the following claims.

\begin{cl}\label{claim1}
If $G$ is residually finite, then we have that either $m\mid n$ or $n \mid m$.
\end{cl}
\begin{proof}
Let $s = \gcd(m,n)$ and let $m',n' \in \mathbb{Z}$ be such that $m = m's$ and $n = n's$. By Bezout's Identity, there exist integers $a$ and $b$ such that $0 \leq a < m'$, $|b| \leq |n'|$ and $an' + 1 = bm'$. 
We consider a group homomorphism $f: G \rightarrow H$ where $H$ is finite. Let $r = \text{ord} (f(x^s))$.

First, we prove that $r$ is co-prime to both $n'$ and $m'$. If there exist $p\mid\gcd(m',r)$, we have 
\begin{equation}\label{eq:ord_co}
1 = f(t^{-1}x^{s\frac{m'r}{p}}t) = f(t^{-1}x^{m\frac{r}{p}}t) =f(x^{s\frac{n'r}{p}}),
\end{equation}
hence we have $p\mid n'$, and since $m'$ and $n'$ are co-prime, we have that $p = 1$. A symmetric argument applies for $\gcd(n',r)$.

Since $\gcd(m',r) = 1$, by Bezout's Identity, there exist $u,v \in \mathbb{Z}$ such that $ur +1 = vm'$. Therefore, we have
\[f(t^{-1}x^st) = f(t^{-1}x^{s(1+ uh)}t) = f(t^{-1}x^{mv}t) =f(x^{nv}) = f(x^s)^{n'v}.\]

Finally, we consider $w = x^s(t^{-1}x^{as}tx^{-bs})^{m'} \in G$. We then have 
\begin{equation}\label{eq:word_neq1}
f(w) = f(x^s)(f(t^{-1}x^st)^af(x^s)^{-b})^{m'} = f(x^s)^{1+an'vm'-bm'} = f(x^s)^{1+an'-bm'}= 1.    
\end{equation}

This shows that for any homomorphism from $G$ to a finite group $H$, the image of $w$ is trivial. Since $G$ is assumed to be residually finite, we have that $w$ must be trivial in $G$. By Britton's Lemma, we have $w = 1$ only if $m'\mid a $ or $n'\mid b$ and so from $an'+1=bm'$, we deduce that that $m' = 1$ or $|n'| = 1$, which is $m\mid n$ or $n \mid m$.
\end{proof}

Since we have that $BS(m,n) \cong BS(n,m)$, without loss of generality, we further assume that $m\mid n$ and $n = mq$ for some $q \in \mathbb{Z}$.

\begin{cl} \label{claim2}
If $G$ is residually finite and $G=BS(m,mq)$, then either $m = 1$ or $q = \pm1$.
\end{cl}
\begin{proof}
For a group homomorphism $f: G \rightarrow H$ with $H$ finite, let $l = \text{ord}(f(t))$. We have $f(x^m) = f(t^{-l}x^mt^l) = f(x^{m})^{q^l}$, so $f(t^{-1}x^{mq^{l-1}}t) = f(x^m)^{q^l} = f(x^m)$ hence we have $f(tx^mt^{-1}) = f(x^{m})^{q^{l-1}}$.

Then we consider $w = x^{-m-q}(tx^mt^{-1}x)^q$. We then have
\begin{equation}\label{eq:word_eq1}
f(w) = f(x)^{-m-q}(f(tx^mt^{-1})f(x))^q = f(x)^{-m-q + q(mq^{l-1}+1)} = f(x)^{mq^l-m} = 1.
\end{equation}

This shows that for any homomorphism from $G$ to a finite group $H$, the image of $w$ is trivial. Since $G$ is assumed to be residually finite, we have that $w$ must be trivial in $G$. By Britton's Lemma, we have $w = 1$ only if $q = \pm 1$ or $m = 1$.
\end{proof}

This shows that if $G$ is residually finite, then either $m=1$ or $|n|=m$.

Conversely, we next prove that if either $m=1$ or $|n|=m$, then we have the conditions $(a)$ or $(b)$ in Lemma \ref{lem:BS_fr}.

\begin{cl}\label{claim3}
There is a normal subgroup $\langle x^k\rangle $ in $G$ for some integer $k\ge 1$ if and only if $|m| = |n|$.
\end{cl}
\begin{proof}
If $|m| = |n|$, it is clear that condition (a) holds as $\langle x^m\rangle$ is a normal subgroup. Conversely, let $k\geq 1$ be an integer such that $\langle x^k\rangle $ is a normal subgroup of $G$. We have $t^{-1}x^{k}t = x^{\varepsilon k}$ where $\varepsilon = \pm 1$ by Lemma \ref{lem:cyclic_normal_subgroups}. Hence we have $t^{-1}x^ktx^{-\varepsilon k} = 1$, by Britton's Lemma we must have $m \mid k$, say $k = mp$. Then, we have $x^{np -\varepsilon k}  = 1$, so $np = \varepsilon k = \varepsilon mp$ and so $|m| = |n|$.
\end{proof}

\begin{cl} \label{claim4}
We have $BS(m,n)$ is solvable if and only if $\min(|m|,|n|) = 1$.
\end{cl}
\begin{proof}
Assume first that $|m| \neq 1$ and $|n|\neq 1$. The subgroup $F$ generated by $a = xtx^{-1}$ and $b = t$ is a free group. Indeed, for every $g\in F$ we can write it in the form 
\[g = a^{n_0}b^{m_1}a^{n_1}b^{m_2}a^{n_2}\cdots b^{m_k}a^{n_k}b^{m_0} = xt^{n_0}x^{-1}(t^{m_1}xt^{n_1}x^{-1})(t^{m_2}xt^{n_2}x^{-1})\cdots(t^{m_k}xt^{n_k}x^{-1})t^{m_0},\] 
where $m_0,n_0 \in \mathbb{Z}$ and $m_i,n_i \in \mathbb{Z}\setminus\{0\}$ for $1\leq i\leq k$ where $k\in \mathbb{N}$. By Britton's lemma, this form is unique, so $F$ is free. Therefore, $BS(m,n)$ contains a free subgroup and so it is not solvable.

Conversely, without loss of generality, we assume that $m = 1$. We then have an exact sequence
\[1 \longrightarrow \llangle x \rrangle \longrightarrow BS(1,n) \longrightarrow \langle t\rangle \longrightarrow 1.\]
Consider the group homomorphism $f: \mathbb{Z}[\frac{1}{n}] = \langle y_0,y_2, \dots \mid y_{i+1}^n = y_i, i \in \mathbb{N} \rangle\rightarrow \llangle x \rrangle$, where for every $k \in \mathbb{N}, f(y_k) = t^kxt^{-k}$. The map $f$ is clearly surjective. We notice that every element $g\in \mathbb{Z}[\frac{1}{n}]$ can be uniquely written as $g = y_0^{n_0}y_1^{n_1}\cdots y_k^{n_k}$, where $n_0 \in \mathbb{Z}, k\in \mathbb{N}, 0 \leq n_i < |n|$ for $1 \leq i \leq k$ and $n_k \neq 0$. With this form, we have that
\[f(g) = x^{n_0}(tx^{n_1}t^{-1})(t^2x^{n_2}t^{-2})\cdots(t^kx^{n_k}t^{-k}) = x^{n_0}tx^{n_1}tx^{n_2}tx^{n_3}\cdots tx^{n_k}t^{-k}.\]
By Britton's Lemma, we have $f(g) = 1$ if and only if $g = 1$, hence $f$ is injective. We conclude that $f$ is an isomorphism. Therefore, $\llangle x \rrangle \cong \mathbb{Z}[\frac{1}{n}]$ and $\langle t\rangle \cong \mathbb{Z}$ are abelian, and so $BS(1,n)$ is metabelian, hence solvable.
\end{proof}

\begin{proof}[Lemma \ref{lem:BS_fr}]
The proof follows from Claims \ref{claim1}, \ref{claim2}, \ref{claim3}, and \ref{claim4}.
\end{proof}

We now prove the converse, namely

\begin{lem}\label{lem:GBS_rf_converse}
Let $G$ be a BS group such that
\begin{enumerate}[(a)]
\item  either for some integer $k\geq 1$, $\langle x^k\rangle $ is a normal subgroup, or
\item $G$ is a solvable BS group, i.e. $G=BS(1,n)$.
\end{enumerate} Then $G$ is linear.
\end{lem}

We first prove that $BS(1,n)$ is linear.
\begin{prop}\label{prop:SolvableBS_linear}
The group homomorphism $f:
 \left\{ \begin{array}{rcl}
 BS(1,n) & \longrightarrow & GL(2,\mathbb{Q}) \\
 x & \mapsto & A \\
 t & \mapsto & T
 \end{array}\right.$, where $A = \left( \begin{array}{cc}
    1 & 1 \\
    0 & 1
 \end{array}\right)$ and $T = \left( \begin{array}{cc}
    \frac{1}{n} & 0 \\
    0 & 1
 \end{array}\right)$ is well defined and injective.
\end{prop}
\begin{proof}
First, we have that $T^{-1}AT = \left( \begin{array}{cc}
    1 & n \\
    0 & 1
 \end{array}\right) = A^n$ so $f$ is well defined.

Next, since $xt = tx^n$ and $t^{-1}x = x^nt^{-1}$, we have for every $k \in \mathbb{Z}, x^kt = tx^{nk}$ and $t^{-1}x^k = x^{nk}t^{-1}$. Therefore, each $g\in BS(1,n)$ can be written in an unique form as $g = t^ix^kt^{-j}$ where $i,j \in \mathbb{N}, k \in \mathbb{Z}$ and if $ij >0$, we have $k\nmid n$ . We then have 
\[f(g) = f(t)^if(x)^kf(t)^{-j} = \left( \begin{array}{cc}
    n^{-i} & 0 \\
    0 & 1
 \end{array}\right)
\left( \begin{array}{cc}
    1 & k \\
    0 & 1
 \end{array}\right)
 \left( \begin{array}{cc}
    n^{j} & 0 \\
    0 & 1
 \end{array}\right)= \left( \begin{array}{cc}
    n^{j-i} & n^{-i}k \\
    0   & 1
 \end{array}\right).\]
So we have $f(g) = 1 \Leftrightarrow i = j = k = 0 \Leftrightarrow g =1$.
So $f$ is injective.
\end{proof}

\begin{proof}[Lemma \ref{lem:GBS_rf_converse}]
The proof is a consequence of Corollary \ref{cor:(a)} and Proposition \ref{prop:SolvableBS_linear}.  
\end{proof}

Lemma \ref{lem:BS_fr} and Lemma \ref{lem:GBS_rf_converse} prove Theorem \ref{thm:main} for Baumslag-Solitar groups.

\section{A single vertex and a finite number of self-loops}
In this section, we suppose that the GBS group $G$ is the fundamental group of a graph of groups with a single vertex and a finite number of self-loops. This group is defined by a family of pairs $(m_i,n_i)_{i\in I}$ of nonzero integers with $I$ finite, which characterize the embeddings of the edge groups into the vertex groups. If $|I| = 1$, then $G$ is a BS group and this case was discussed in the previous section. We further suppose that $|I|\geq 2$.

\begin{prop} \label{prop: mult. loops}
$G$ is residually finite if and only if $|n_i|=|m_i|$ for all $i\in I$.   
\end{prop}

\begin{proof}

We first assume that $G$ is residually finite. Since a subgroup of a residually finite group is residually finite, using the case $|I| = 1$ proven in the previous section, we have that for every $i\in I$, either $|m_i|=|n_i|$ or $\min(|m_i|,|n_i|)=1$. We assume without loss of generality that whenever $\min(|m_i|,|n_i|)=1$, we have $n_i = 1$ or $m_i = 1$.

\begin{cl} \label{claim5}
For each $i\in I$, $(m_i,n_i)$ is not $(1,n)$ or $(m,1)$, where $|m|,|n|\geq 2$.
\end{cl}
\begin{proof}
Suppose that there exists such $(1,n)$ $|n|\geq 2$, among the $(m_i,n_i)_{i \in I}$. We consider the three possible cases:\\
Case 1: There is another pair $(1,m)$ among the pairs $(m_i, n_i)_{i\in I}$. Then $G$ has a subgroup with presentation 
\[\langle t,u,x \mid t^{-1}xt=x^n,u^{-1}xu=x^m\rangle.\]
For a group homomorphism $f: G \rightarrow H$ with $H$ finite, let $r = \text{ord} (f(x))$. We must have that $\gcd(r,m) = \gcd(r,n) = 1$ as we have seen in equation \eqref{eq:ord_co}. By Bezout's Identity, there exist $a,b \in \mathbb{Z}$ such that $an \equiv bm \equiv 1 (\text{mod }r)$. We have that the commutator $ [txt^{-1},uxu^{-1}]$ is not trivial in $G$ by Britton's Lemma. However,
\[f([txt^{-1},uxu^{-1}]) = [f(txt^{-1}),f(uxu^{-1})] = [f(tx^{an}t^{-1}), f(ux^{bm}u^{-1})] = [f(x^{a}),f(x^{b})] = 1.\]
Hence, in this case $G$ is not residually finite.\\
Case 2: There is a pair $(m,1)$ among the other pairs. Then $G$ has a subgroup with presentation 
\[\langle t,u,x \mid t^{-1}xt=x^n,uxu^{-1}=x^m\rangle.\]
This subgroup is isomorphic to $\langle t,u,x \mid t^{-1}xt=x^n,u^{-1}xu=x^m\rangle$, which is Case 1. Hence $G$ is not residually finite.\\
Case 3: The other pairs are of the form $|n_i|=|m_i|$. We consider one of the form $(m, em)$ with $e = \pm 1$. Then $G$ has a subgroup with presentation 
\[\langle t,u,x \mid t^{-1}xt=x^n,u^{-1}x^mu=x^{em}\rangle.\]
For a group homomorphism $f: G \rightarrow H$ with $H$ finite, let $r = \text{ord} (f(x))$. Similarly, there exist $a \in \mathbb{Z}$ such that $an \equiv 1 (\text{mod }r)$.
We see that $utx^{m}t^{-1}u^{-1}tx^{-em}t^{-1}$ is not trivial in $G$ by Britton's Lemma, as if it is, we must have $n\mid m$ and then $1 = utx^{m}t^{-1}u^{-1}tx^{-em}t^{-1} = ux^{\frac{m}{n}}u^{-1}x^{-e\frac{m}{n}}$, and again with Britton's Lemma, it implies that $n = 1$ contradicting the assumption that $|n|\geq 2$. Then we have that
\begin{equation} \label{eq:word_2loops}
\begin{aligned}
f(utx^{m}t^{-1}u^{-1}tx^{-em}t^{-1}) &= f(u)f(tx^{m}t^{-1})f(u)^{-1}f(tx^{-em}t^{-1})\\
&= f(u)f(tx^{anm}t^{-1})f(u)^{-1}f(tx^{-eanm}t^{-1}) \\
&= f(u)f(x^{am})f(u)^{-1}f(x^{-eam})\\
&= f(ux^{am}u^{-1}x^{-eam}) = 1.
\end{aligned}
\end{equation}
Hence, $G$ is not residually finite.
\end{proof}
Conversely, suppose that every pair $(m_i, n_i)_{i\in I}$ satisfies that $|m_i| = |n_i|$. Let $k = \text{lcm}(m_i)$. Then, $\langle x^k\rangle$ is contained in all edge groups, and hence it is a normal subgroup of $G$. 
\end{proof}

This proves Theorem \ref{thm:main} for the case of a single vertex and a finite number of self-loops.

\section{General GBS groups}
For the general cases, we start by introducing a property of graph of groups, which we prove characterizes linearity for a GBS group.
\begin{mydef}[Unimodular] \label{def:unimodular}
Given a sequence of oriented edges $e_1,\dots,e_p$ forming a loop ($p\geq 1$). Let the left and right inclusion of $e_j$ be given by multiplication by $n_j, m_j$, respectively. We say that the loop is \emph{unimodular} if $|\prod^{p}_{j=1} m_j| = |\prod^{p}_{j=1} n_j|$.\\
We say that the graph of groups of a GBS group is unimodular if all loops are unimodular. In this case, we say that the group is unimodular.
\end{mydef}
This property can be easily read through the graph of groups.

\begin{rem} \label{rem:subgroups}
Let $G$ be a GBS group. Consider a subgroup which is the fundamental group of a simple loop $\ell$ in the graph of groups, defined by the edges $e_1,\dots,e_p$($p\geq 1$). Let the left and right inclusion of $e_j$ be given by multiplication by $n_j, m_j$, then $G$ has a subgroup with presentation
\[H_\ell=\langle t_1,\dots, t_p, x_1,\dots, x_p\mid t_i =1,i\in I, t_k^{-1}x_k^{m_k}t_k = x_{k+1}^{n_k},1\leq k<p, t_p^{-1}x_p^{m_p}t_p = x_1^{n_p}\rangle,\]
where $I$ depends on the spanning tree. Furthermore, the subgroup generated by $u = t_1t_2\dots t_p$ and $x_1$, by Britton's lemma, satisfies
\[u^{-1}x_1^{\prod^{p}_{j=1} m_j}u = x_1^{\prod^{p}_{j=1} n_j},\]
and is a Baumslag-Solitar group with presentation
\[BS_\ell=\langle x_1,u \mid u^{-1}x_1^{\prod^{p}_{j=1} m_j}u = x_1^{\prod^{p}_{j=1} n_j} \rangle.\]    
\end{rem}

\begin{lem} \label{lem:unimodular}
The GBS group $G$ is unimodular if and only if for some integer $k\geq 1$, $\langle x^k\rangle $ is a normal subgroup.
\end{lem}
\begin{proof}

From Remark \ref{rem:subgroups}, we have that for each loop $\ell$ in the graph of groups, we can consider the corresponding Baumslag-Solitar group $BS_\ell$. From Claim 3, we have that if $G$ has a normal subgroup $\langle x^k\rangle $ for some $k\geq 1$, we must have that the loop is unimodular. This should hold for each loop, so $G$ is unimodular.

Conversely, suppose that $G$ is unimodular and let $k\in\mathbb{N}$ be the product of each integer $n_j$, $m_j$ defined for the inclusion of the edge groups. Then $x^k$ is an element in all vertex groups. For each generator $t_e$ associated to an edge $e$, we have that either $t_e = 1$ (that is, the edge $e$ is contained in the spanning tree), or there exists a loop $e = e_1,e_2,\dots,e_p$, where all edges other than $e_1$ are contained in the spanning tree. In both cases, we have that $t_e$ satisfies $t_e^{-1}x^kt_e = x^{\pm k}$ and so it follows that $\langle x^k \rangle$ is a normal subgroup of $G$.
\end{proof}

\begin{lem} \label{lem:rf_cases}
If $G$ is residually finite, then one of the following holds
\begin{enumerate}[(a)]
\item  either for some integer $k\geq 1$, $\langle x^k\rangle $ is a normal subgroup, or
\item $G$ is a solvable BS group, i.e. $G=BS(1,n)$.
\end{enumerate}
\end{lem}
\begin{proof}

Consider the subgroup which is the fundamental group of a loop $\ell$ as in Remark \ref{rem:subgroups} and the associated $BS_\ell$. From Claim \ref{claim3}, we have that $BS_\ell$ is either solvable or modular, and hence the loop $\ell$ is unimodular, or the products of one of the two multiplications (right or left) must be $\pm 1$.

For the case where there are more than $2$ loops, it follows from Section 4 that the group $G$ is unimodular; hence, by Lemma \ref{lem:unimodular}, condition (a) holds.

Suppose next that there is only one loop in the graph. If the loop is unimodular, then we conclude in virtue of Lemma \ref{lem:unimodular}. We suppose that the loop is not unimodular. Without loss of generality, we suppose that $m_k = 1$ for every $k$ and its fundamental group is
\[\langle x_1,x_2, \dots, x_k,t \mid x_i = x_{i+1}^{n_i}, 1 \leq i < k, t^{-1}x_kt = x_1^{n_k}\rangle.\]
This group is clearly just 
\[BS(1,\prod_{j=1}^kn_j) = \langle x_k, t\mid t^{-1}x_kt = x_k^{\prod_{j=1}^kn_j}\rangle.\]
Hence, we can suppose that the single loop is a self-loop. If there are other vertices in the graph, we are in the case (b). Otherwise, $G$ has a subgroup with presentation (for some $l,m,n$)
\[\langle x,y,t\mid t^{-1}xt = x^l, x^m =y^n\rangle.\]

\begin{cl}\label{prop:subgroup_loopedge}
The subgroup $K$ generated by $y,t \in H = \langle x,y,t\mid t^{-1}xt = x^l, x^m =y^n\rangle$ has a presentation $\langle y,t\mid t^{-1}y^nt = y^{nl}\rangle$.
\end{cl}
\begin{proof}
We consider an HNN extension $H$ of the group $\langle x,y\mid x^m =y^n\rangle$ relative to the isomorphism $\varphi:
 \left\{ \begin{array}{rcl}
 \langle x\rangle & \longrightarrow & \langle x^l\rangle \\
 x & \mapsto & x^l
 \end{array}\right.$. Now, the subgroup $K$ acts on the Bass-Serre covering tree of $H$ and so it is the fundamental group of a graph of cyclic vertex and edge groups. Since $t$ is an edge element, this graph of groups has one vertex and a self-loop. We have that $K = \langle y,t \mid t^{-1}y^at = y^b\rangle \cong BS(a,b)$ for some $a,b \in \mathbb{Z} \setminus \{0\}$. Since we have $t^{-1}y^nt =t^{-1}x^mt = x^{ml} =y^{nl}$ in $H$, we have $a \mid n$ and $b = al$, by Britton's lemma we must have $a=n$.
\end{proof}

Since $G$ is assumed to be residually finite, so is the BS group $K$, and so we have that either $n = 1$ or $l = 1$. 

Since we are in the case that $G$ is not unimodular, we have that $l$ is not $1$. Therefore, we have that $n = 1$, and in this case we have that $y \in \langle x\rangle$, so every edge outside the loop can be embedded in the vertex with the self loop, hence $H$ is isomorphic to $BS(1,l)$.
\end{proof}

\begin{proof}[Theorem \ref{thm:main}]
The proof follow from Corollary \ref{cor:(a)}, Lemma \ref{lem:BS_fr}, and Lemma \ref{lem:rf_cases}.
\end{proof}
\section{Generalization: Infinite virtually cyclic groups}

In this section, we consider the fundamental group of a graph of groups with infinite virtually cyclic vertex and edge groups. 

We start with a few observations on infinite cyclic groups.

\begin{lem} \label{lem:infcyc_intersect}
For an infinite virtually cyclic group $G$, if $g \in G$ is such that $\langle g\rangle$ is a finite index subgroup of $G$, then every infinite cyclic subgroup $\langle h\rangle$ of $G$ satisfies $\langle g\rangle\cap  \langle h \rangle \neq \{1\}$.
\end{lem}
\begin{proof}
We consider a partition with $n = [G:\langle g\rangle]$ and $g_1 = 1,g_2,\dots,g_n \in G$:
\[G = \bigsqcup_{i = 1}^{n}g_i\langle g \rangle.\]
Consider the left cosets $\langle g\rangle, h\langle g\rangle,h^2\langle g\rangle, \dots h^n\langle g\rangle$. By the pigeonhole principle, we have that at least two of the cosets are equal, and hence we have $\langle g\rangle\cap \langle h \rangle \neq \{1\}$.
\end{proof}

By the multiplicity of index, we have the following
\begin{cor} \label{cor:infcyc_finite_index}
Every infinite cyclic subgroup of an infinite virtually cyclic group has finite index.
\end{cor}

\begin{lem} \label{lem:infcyc_normal}
For an infinite virtually cyclic group $G$, if $g \in G$ is such that $\langle g\rangle$ is a finite index subgroup of $G$, then there exists $n>0$ such that $\langle g^n\rangle \trianglelefteq G$.
\end{lem}
\begin{proof}
We consider a partition with $m = [G:\langle g\rangle]$ and $g_1 = 1,g_2,\dots,g_m \in G$:
\[G = \bigsqcup_{i = 1}^{m}g_i\langle g \rangle.\] Consider the intersection
\[\langle g^n\rangle = \bigcap_{i = 1}^mg_i\langle g\rangle g_i^{-1}.\]
Notice that $n>0$ exists by Lemma \ref{lem:infcyc_intersect} and this subgroup is clearly normal.
\end{proof}

\begin{lem} \label{lem:infcyc_int_vc}
Let $G$ be an infinite virtually cyclic group and let $H$ be an infinite virtually cyclic subgroup of $G$. If $g \in G$ is such that $\langle g\rangle$ is an infinite subgroup of $G$, then we have $\langle g^m \rangle\subseteq H$ for some $m > 0$.
\end{lem}
\begin{proof}
By Lemma \ref{lem:infcyc_normal}, we have an exact sequence
\[1 \rightarrow \langle g^n\rangle \rightarrow G \xrightarrow{\pi} Q \rightarrow 1,\]
where $Q$ is a finite group. Then we consider the intersection with $H$; we have another exact sequence
\[1 \rightarrow \langle g^n\rangle \cap H \rightarrow H \xrightarrow{\pi} \pi(H) \rightarrow 1.\]
If $\langle g\rangle \cap H = \{1\}$, then we have $\langle g^n\rangle\cap H = \{1\}$ and $H \cong \pi(H) \subseteq Q$ is finite, contradiction. So $\langle g\rangle \cap H \neq \{1\}$ hence $\langle g^m \rangle\subseteq H$ for some $m > 0$.
\end{proof}

With these elements, we prove the generalized case of the fundamental group of a graph of groups with infinite virtually cyclic vertex and edge groups.

\begin{proof}[Theorem \ref{thm:v.c.}]
Let $T$ be a spanning tree of the graph of groups and $G_T\subseteq G$ the corresponding fundamental group. By Lemma \ref{lem:infcyc_normal} and Lemma \ref{lem:infcyc_int_vc}, $G_T$ clearly verifies (a). If there are no additional edges outside of $T$, then $G = G_T$ verifies (a).

We start by proving the fact that if $G$ is residually finite, then we have either condition (a) or (b).

Suppose that there is an additional edge $e$ equipped with the group isomorphism $\varphi:H \rightarrow K$ and the edge element $t$. We then have $G$ has a subgroup
\[G_T*_{\varphi} =\langle G_T,t \mid t^{-1}ht = \varphi(h), h\in H\rangle.\]

For every $g'\in G_T$ of infinite order contained in a vertex group, we have $\langle g' \rangle \cap \langle tg't^{-1} \rangle \neq \{1\}$. We consider $m'> 0$ such that ${g'}^{m'}$ generates $\langle g' \rangle \cap \langle tg't^{-1} \rangle$. So, there exists $n' \in \mathbb{Z}_{\neq 0}$ such that $t^{-1}{g'}^{m'}t = {g'}^{n'}$. Taking $d = \gcd(m',n'), m = \frac{m'}{d}, n= \frac{n'}{d},$ and $g_0 = {g'}^d$, we have $t^{-1}g_0^mt = g_0^n$ and $\gcd(m,n) = 1$. We wish to prove that $m = 1$ or $|n| = 1$, we suppose towards contradiction that $m > 1, |n| > 1$.

We construct a sequence $(g_i)_{i\geq 0}$ in $G_T$ such that $t^{-1}g_i^mt = g_i^n$ and $[H:\langle g_{i+1}^m \rangle] < [H:\langle g_{i}^m \rangle]$ by induction, with $g_0$ given in the previous paragraph. Suppose that $g = g_i$ is constructed; we construct $g_{i+1}$ following the idea and calculation in Claim \ref{claim1}. Let $a, b \in \mathbb{Z}$ verify Bezout's Identity $an + 1 = bm$. Consider a group homomorphism $f:G_T*_{\varphi}\rightarrow F$ with $F$ finite and let $r = \text{ord}(f(g))$.

First, we have that $r$ is co-prime to both $m$ and $n$ from the same calculation with equation \eqref{eq:ord_co}. By Bezout's Identity, there exists $u_1,v_1,u_2,v_2 \in \mathbb{Z}$ such that $u_1r + 1 = v_1m$ and $u_2r + 1 = v_2n$. Therefore, we have
\[f(t^{-1}gt) = f(g)^{nv_1} \text{ and } f(tgt^{-1}) = f(g)^{mv_2}.\]
Then consider $w =  g(t^{-1}g^{a}tg^{-b})^{m} \in G_T*_\varphi$, we have $f(w) = 1$ by equation \eqref{eq:word_neq1}, so by residually finiteness of $G_T*_\varphi$, we have that $w = 1$. By Britton's lemma, we have two possible cases:

Case 1: $g^a \in H$, with $g^m \in H$ and $\gcd(a,m)=1$ we have $g \in H$. Let $k = t^{-1}gt \in K$, we have $f(k) = f(g)^{nv_1}$ and $g(k^ag^{-b})^{m} = 1$. We define $g_{i+1} = g^bk^{-a} \in G_T$, we have $g_{i+1}^m = g = g_i$ so $[H:\langle g_{i+1}^m \rangle] = \frac{1}{m}[H:\langle g_{i}^m \rangle] < [H:\langle g_{i}^m \rangle]$. It remains to prove that $t^{-1}g_{i+1}^mt = g_{i+1}^n$, notice that
\[f(kg_{i+1}^{-n}) = f(k)(f(k)^af(g)^{-b})^n = f(g)^{an^2v_1 + nv_1 - bn} = f(g)^{bmnv_1 - bn} =f(g)^{bnu_1h} = 1,\]
by residually finiteness we have $g_{i+1}^n = k = t^{-1}gt = t^{-1}g_{i+1}^mt$.

Case 2: $g^a \notin H$, which implies that $g^b \in K$, with $g^n \in K$ and $\gcd(b,n)=1$ we have $g \in K$. Let $h = tgt^{-1} \in H$, we have $f(h) = f(g)^{mv_2}$ and
\[1 = g(t^{-1}g^{a}tg^{-b})^{m}= tg(t^{-1}g^{a}tg^{-b})^{m}t^{-1} = h(g^ah^{-b})^m.\]
Let $g_{i+1} = h^{b}g^{-a} \in G_T$, notice that
\[f(gg_{i+1}^{-n}) = f(g)(f(g)^af(h)^{-b})^n = f(g)^{an + 1 - bnmv_2} = f(g)^{an+1 -bm} = 1,\]
by residually finiteness we have $g_{i+1}^n = g$, with $g_{i+1}^m = h \in H$, we have $[H:\langle g_{i+1}^m \rangle]  < [H:\langle g_{i}^m \rangle]$. It remains to prove that $t^{-1}g_{i+1}^mt = g_{i+1}^n$, which is clear since $tg_{i+1}^nt^{-1} = tgt^{-1} = h = g_{i+1}^m$.

In both cases, we construct $g_{i+1}$ to verify the condition, so we conclude that such $(g_i)_{\geq 0}$ exists by induction. This leads to a contradiction since $[H:\langle g_{i}^m \rangle] \in \mathbb{Z}_{>0}$ for all $i$. So we must have $m = 1$ or $|n| = 1$.

This leads to four cases:\\
Case 1: $\varphi(g^m) = g^{mq}$, for some $m \in \mathbb{Z}_{>0},  q \in \mathbb{Z}\setminus\{-1,0,1\}.$\\
Case 2: $\varphi^{-1}(g^{m}) = g^{mq}$, for some $m \in \mathbb{Z}_{>0}, q \in \mathbb{Z}\setminus\{-1,0,1\}.$\\
Case 3: $\varphi(g^m) = g^m$, for some $n \in \mathbb{Z}\setminus \{0\}$\\
Case 4: $\varphi(g^m) = g^{-m}$, for some $n \in \mathbb{Z}\setminus\{0\}$.

By Lemma \ref{lem:infcyc_intersect}, for each $g,h \in G_T$ of infinite order contained in a vertex group, the corresponding cyclic subgroups $\langle g\rangle$ and $\langle h\rangle$ intersect nontrivially, and so both $g$ and $h$ satisfy the same alternative above and, in cases 1 and 2, they satisfy it with the same $q$. It follows that each element of infinite order contained in a vertex group satisfies the same alternative. Notice that cases 1 and 2 are the same since we can replace $\varphi$ with $\varphi^{-1}$, which will not change the structure of the group.

Case 1(and 2): Suppose that there exists a fixed $q \in \mathbb{Z}\setminus\{-1,0,1\}$ such that for each $g \in G_T$ of infinite order contained in a vertex group, we have $\varphi(g^m) = g^{mq}$ for some $m \in \mathbb{Z}_{>0}$. In this case, we prove that $H = G_T$. We suppose, towards contradiction, that $G_T \neq H$. Let $g\in G_T \setminus H$ of infinite order and contained in a vertex group. Let $m \in \mathbb{Z}_{>0}$ be such that $\varphi(g^m) = g^{mq}$. Then, define the sequence $(h_l)_{l \in \mathbb{Z}_{\geq 0}}$ in $G_T*_{\varphi}$ by $h_0 = g^m$ and $h_{l+1} = th_lt^{-1}$. If $h_l \in H$, then since $h_l = t^lg^mt^{-l}$, we have $h_l^{q^l} = t^lg^{mq^l}t^{-l} = g^m$, so that $[H:\langle h_l\rangle] = \frac{[H : \langle h_l^{q^l}\rangle]}{|q|^l} = \frac{[H:\langle g^m\rangle]}{|q|^l}$. Since the index must be an integer, with $|q| > 1$, there exists $n\in \mathbb{Z}_{\geq 0}$ such that $h_{n} \in H$ and $h_{n+1} \notin H$, hence $h_n \in G_T\setminus K$. We consider $w = g^{-m-q^{n+1}} (th_nt^{-1}g)^{q^{n+1}}=g^{-m-q^{n+1}} (t^{n+1}g^mt^{-n-1}g)^{q^{n+1}}$, with similar calculation as in equation \eqref{eq:word_eq1}, we have that for any homomorphism from $G_T*_\varphi$ to a finite group, the image of $w$ is trivial, by residually finiteness $w = 1$, which contradicts Britton's Lemma, so we must have that each element of $G_T$ of infinite order and contained in a vertex group must belong to $H$.

Let $g_0\in G_T$ be of infinite order contained in a vertex group and such that $\langle g_0\rangle \trianglelefteq G_T$, which exists since $G_T$ verifies the alternative (a). For $a \in G_T$, by Lemma \ref{lem:cyclic_normal_subgroups} we have that $ag_0a^{-1} = g_0^e$ with $e = \pm 1$, and by Corollary \ref{cor:infcyc_finite_index}, we have $\langle g_0\rangle\subseteq H$ is finite index. We suppose that $g \in H$ is the element such that $\langle g\rangle$ is of minimal index such that $aga^{-1} = g^e$ and $\varphi(g) = g^q$. Then we consider $atgt^{-1}a^{-1}tg^{-e}t^{-1} \in G_T*_{\varphi}$, and with the same calculation as in equation \eqref{eq:word_2loops}, we have that $f(atgt^{-1}a^{-1}tg^{-e}t^{-1}) = 1$ for a group homomorphism $f:G_T*_{\varphi}\rightarrow F$ with $F$ finite, so we must have $atgt^{-1}a^{-1}tg^{-e}t^{-1} = 1$. If $a\notin H$, by Britton's lemma, we must have that $g\in K$, hence $\varphi(h) =g$ for some $h\in H$. Since $\varphi(h^q) = g^q$, we have $\varphi(h) = g = \varphi^{-1}(g^q)= h^q$ and $aha^{-1} =atgt^{-1}a^{-1} = tg^et^{-1} = h^e$  with $[H:\langle h\rangle] = \frac{[H:\langle g\rangle]}{[\langle h\rangle: \langle g\rangle]} < [H:\langle g\rangle]$, which contradicts the minimality of the index of the subgroup generated by $g$. Hence $a \in H$. We conclude that $H = G_T$ and so we have that $G_T$ is virtually cyclic. Therefore, we can suppose that $T$ is a single vertex with $H$ being its vertex group.

We prove that this is the only self-loop. We suppose, towards contradiction, that there exists another edge with edge element $u$ and corresponding group homomorphism $\psi: H' \rightarrow K'$. Then $\psi$ also falls into one of the four alternative cases above.

Subcase 1(and 2): Suppose that there exists a fixed $p \in \mathbb{Z}\setminus\{-1,0,1\}$ such that for each $g \in H$ of infinite order contained in a vertex group, we have $\psi(g^m) = g^{mp}$ for a $m \in \mathbb{Z}_{>0}$. We have already proved that in this case $H = H'$. Notice that since there exists $g\in H$ such that $u^{-1}gu = g^p$ and $[H:\langle g\rangle] = [\psi(H): \langle \psi(g) \rangle]= [K':\langle g^p\rangle] < +\infty$, we have that $K'\neq H$ and there exist $h \in H \setminus K'$ of infinite order and such that $[H: \langle h\rangle]$ is minimal. Then $k = uhu^{-1} \notin H$. Notice that there exists $m > 0$ such that $t^{-1}k^mt = k^{mq}$. Then, define the sequence $(h_l)_{l \in \mathbb{Z}_{\geq 0}}$ in $G$ by $h_0 = k^m$ and $h_{l+1} = th_lt^{-1}$. If $h_l \in H$, then since $h_l = t^lk^mt^{-l}$, we have $h_l^{q^l} = t^lk^{mq^l}t^{-l} = k^m$, so that $[H:\langle h_l\rangle] = \frac{[H : \langle h_l^{q^l}\rangle]}{|q|^l} = \frac{[H:\langle k^m\rangle]}{|q|^l}.$ Since the index must be an integer, with $|q| > 1$, there exists $n\in \mathbb{Z}_{\geq 0}$ such that $h_{n} \in H$ and $h_{n+1} \notin H$, hence $h_n \in H\setminus K$. We consider $w = k^{-m-q^{n+1}} (tkt^{-1}k)^{q^{n+1}}=k^{-m-q^{n+1}} (t^{n+1}k^mt^{-n-1}k)^{q^{n+1}}$, with similar calculation as in Equation \eqref{eq:word_eq1}, we have that for any homomorphism from $G$ to a finite group, the image of $w$ is trivial, and so by residually finiteness we have that $w = 1$, which contradicts Britton's Lemma.

Subcase 3 and 4: Suppose that there exists $e = \pm 1$ such that for each $g \in H$ of infinite order, we have $\psi(g^m) = g^{em}$ for a $m \in \mathbb{Z}_{>0}$. There exists $g \in H$ such that $t^{-1}gt = g^q$ and $u^{-1}gu = g^{e}$, we suppose that $g\in H$ with minimal $[H:\langle g\rangle]$ among the elements satisfying these two equations. For a group homomorphism $f: G \rightarrow F$ with $F$ finite, let $r = \text{ord} (f(g))$. Similarly proven, there exists $a \in \mathbb{Z}$ such that $aq \equiv 1 (\text{mod }r)$.
We consider $w = utgt^{-1}u^{-1}tg^{-e}t^{-1} \in G$, with the same calculation as in equation \eqref{eq:word_2loops}, we have that $f(w) = 1$, hence $w = 1$ by residually finiteness of $G$. By Britton's Lemma, we have that $h = tgt^{-1} \in H$ and then $1 = uhu^{-1}h^{-e}, u^{-1}hu = h^e$, and $h^q = tg^qt^{-1} = g = t^{-1}ht$, but $[H:\langle g\rangle] = [H:\langle h^q\rangle] = |q|[H:\langle h\rangle] > [H:\langle h\rangle]$, which contradicts the minimality of index in the choice of $h$.

The discussion above implies that there can only be a single self-loop. In this case, we have that $G = H*_\varphi$ is an ascending HNN-extension, which is case (b) in the statement.

Conversely, suppose that $G = H*_\varphi$ is an ascending HNN-extension of an infinite virtually cyclic group relative to $\varphi: H \rightarrow K$. For every $g \in H$ of infinite order, there exists $m,n \in \mathbb{Z}_{\neq 0}$ such that $\varphi(g^m) = g^n$, we have that $[H:K] = \frac{[H:\langle g^{mn}\rangle]}{[K:\langle g^{mn}\rangle]} = \frac{[K:\langle g^{n^2}\rangle]}{[K:\langle g^{mn}\rangle]} = |\frac{n}{m}|$, so we have $m \mid n$. So for every $g \in H$ of infinite order, we have $\varphi(g^m) =g^{mq}$ for some $m \in \mathbb{Z}_{>0},  q \in \mathbb{Z}\setminus\{-1,0,1\}$. Let $g \in H$ of infinite order and $[H:\langle g\rangle]$ minimal such that $\varphi(g) = g^q$. Consider the subgroup $S$ of $H*_{\varphi}$ generated by $g$ and $t$, we wish to prove that it is isomorphic to $BS(1,q) = \langle h,s\mid s^{-1}hs = h^q\rangle$, we construct $f: BS(1,q) \rightarrow S$ by $f(s) = t$ and $f(h) = g$ it is clearly surjective, we prove that it is also injective. Since each element in $BS(1,q)$ can be uniquely written as $s^ih^ks^{-j}$ with $i,j \in \mathbb{N}, k \in \mathbb{Z}$ and if $ij > 0$, we have $k\nmid q$, we have $f(s^ih^ks^{-j}) = 1 \Leftrightarrow t^ig^kt^{-j}=1$, by iterating Britton's lemma we have $g^k \in \varphi^i(H)$ and in this case we have $g^k = 1$ hence $k = 0$ and $s^ih^ks^{-j}$, we conclude $S\cong BS(1,q)$. Hence $S$ is linear by Theorem \ref{thm:main}.

Also, since for all elements $a \in H$ we have $t^{-1}at = \varphi(a)$, every element in $H*_{\varphi}$ can be written uniquely in the form $t^ist^{-j}$ with $s \in H$ and $i,j \in \mathbb{N}$ and if $ij > 0, s\notin K$. We consider the partition of $H = \bigsqcup_{l=1}^{m}\langle g\rangle g_l$ where $g_l \in H, l =1,\dots,m$. We wish to prove that if $l \neq l'$, $\langle g \rangle \varphi(g_l) \neq \langle g \rangle \varphi(g_{l'})$. Suppose towards contradiction that there exist $g_l \notin \langle g \rangle$ such that $\varphi(g_l) = g^n$, then there exists an element $s \in H$ such that $\varphi(s) = g^d$ where $d = \gcd(n, q)$. Let $p =\frac{q}{d}$. Since $\varphi(s^p) = g^{pd} = g^{q}$, we have $s^p = g$ and $\varphi(s) = s^q$, with $[H:\langle s \rangle] = \frac{[H:\langle g\rangle]}{|p|}$, we must have $|p| = 1$, which is $q \mid n$, which contradicts that $g_l \notin \langle g \rangle$. Hence, we have that $H = \bigsqcup_{i=l}^{m}\langle g\rangle \varphi^j(g_l)$ for every $j \in \mathbb{N}$. Since each element in $H*_{\varphi}$ can be uniquely written as $t^ig^k\varphi^j(g_l)t^{-j} = t^ig^kt^{-j}g_l$, we have a partition $H*_{\varphi} = \bigsqcup_{i=1}^{m}Sg_i$. This implies that $S\subseteq G$ is of finite index. Since $S$ is linear, we have that $H*_{\varphi}$ is virtually linear, and by Lemma \ref{lem:v.l.}, we conclude that $G = H*_{\varphi}$ is linear. 

Case 3 and 4: Case 1(and 2) implies that if there is more than one loop, then we must have that for every edge element $t$, there exists $m\in \mathbb{Z}\setminus \{ 0 \}$ satisfying that $t^{-1}g^mt = g^{\pm m}$ for each $g$ in the vertex group of infinite order. Then we can consider $g$ of infinite order such that $\langle g\rangle$ is a normal subgroup of every vertex group and edge group. Then $\langle g\rangle$ is a normal subgroup of $G$. Conversely, an analogous proof to that of Lemma \ref{lem:BS_normal_linear} shows that $G*_{\varphi}$ is virtually a direct product of $\mathbb{Z}$ and a free group, hence virtually linear and, so by Lemma \ref{lem:v.l.}, linear.
\end{proof}

\section*{Acknowledgments}
I would like to thank Montserrat Casals-Ruiz for supervising me during my internship at the University of the Basque Country(UPV/EHU), she helps a lot during the processes of writing, without her help it is impossible for me to finish it. 
\printbibliography[
heading=bibintoc,
title={Bibliography}
]  

@book{Serre,
    author = {Serre, Jean-Pierre},
    title = {Arbres, amalgames, $SL_2$},
    series = {Ast\'erisque},     
    publisher = {Soci\'et\'e math\'ematique de France},
    number = {46},
    year = {1977},
    mrnumber = {476875},
    zbl = {0369.20013},
    language = {fr},
    url = {https://www.numdam.org/item/AST_1983__46__1_0/}
}

@online{Cornulier,
    title = {When is a generalised Baumslag-Solitar group linear?},
    author = {Yves de Cornulier},
    note = {(version: 2022-08-16)},
    howpublished = {MathOverflow},
    url = {https://mathoverflow.net/q/428598}
}

@book{Roger,
    author = {Roger C. Lyndon and Paul E. Schupp},
    title = {Combinatorial Group Theory},
    publisher = {Springer-Verlag},
    year = {2001}
}

@article{Levitt,
    url = {https://doi.org/10.1515/jgth-2014-0028},
    title = {Quotients and subgroups of Baumslag–Solitar groups},
    author = {Gilbert Levitt},
    pages = {1--43},
    volume = {18},
    number = {1},
    journal = {Journal of Group Theory},
    doi = {doi:10.1515/jgth-2014-0028},
    year = {2015},
}
\end{document}